\documentclass[11pt]{article}
\usepackage{amsmath} 
\usepackage{geometry}                
\usepackage[parfill]{parskip}    
\usepackage{graphicx}
\usepackage[noend]{algorithm2e}
\usepackage{amssymb}
\usepackage{epstopdf}
\usepackage{verbatim}
\usepackage{hyperref}
\hypersetup{colorlinks,%
citecolor=black,%
filecolor=black,%
linkcolor=black,%
urlcolor=black%
}
\usepackage{amsthm}
\newtheorem{thm}{Theorem}
\newtheorem{cor}[thm]{Corollary}
\newtheorem{lem}[thm]{Lemma}

\theoremstyle{definition}
\newtheorem{defn}{Definition}

\DeclareMathOperator{\Tr}{Tr}
\DeclareMathOperator{\wt}{wt}

\newcommand{\f}[1]{\mathbb{F}_{#1}}
\newcommand{\tfr}[1]{{\widehat{#1}}}
\newcommand{\K}{\mathcal{K}}

\title{Ternary Kloosterman sums using Stickelberger's theorem and the Gross-Koblitz formula}

\author{Faruk G\"olo\u glu\thanks{Research supported by Claude Shannon Institute,
Science Foundation Ireland Grant 06/MI/006}, Gary McGuire\footnotemark[1]
, Richard Moloney\thanks{Research supported by Claude Shannon Institute,
Science Foundation Ireland Grant 06/MI/006, 
and the Irish Research Council
for Science, Engineering and Technology}\\
School of Mathematical Sciences\\
University College Dublin\\
Ireland\\} 

\begin{document}
\maketitle

\begin{abstract}
We give  results characterising ternary Kloosterman sums modulo $9$ and $27$. This leads to a complete characterisation of values that ternary Kloosterman sums assume modulo $18$ and $54$. The proofs uses Stickelberger's theorem, the Gross-Koblitz formula and Fourier analysis. 
\bigskip

\end{abstract}

\section{Introduction}
Let $\mathcal{K}_{p^n}(a)$ denote the $p$-ary Kloosterman sum defined by
\[
\mathcal{K}_{p^n}(a) := \sum_{x \in \f{p^n}} \zeta^{\Tr (x^{p^n-2}+ax)},
\]
for any $a \in \f{p^n}$, where $\zeta$ is a primitive $p$-th root of unity 
and $\Tr$ denotes the absolute trace map $\Tr : \f{p^n} \to \f{p}$ defined as usual as
\[
\Tr (c) := c + c^p + c^{p^2} + \cdots + c^{p^{n-1}}.
\]

Kloosterman sums have
attracted attention thanks to their various links to other related fields. For instance, 
a zero of a binary Kloosterman sum on $\f{2^n}$ leads to a bent function from $\f{2^{2n}} \to \f{2}$ as proven by Dillon in \cite{Dil74}. Similarly, zeros of ternary Kloosterman sums give rise to ternary bent functions \cite{HK}. However determining a zero of a Kloosterman sum is not easy. A recent result in this direction is the following: a binary or ternary Kloosterman sum $\mathcal{K}_{p^n} (a)$ is not zero if $a$ is in a proper subfield of $\f{p^n}$ 
except when $p = 2, n = 4, a = 1$, see \cite{Lisonek2009}. 
Given the difficulty of the problem of finding zeros (or explicit values) of Kloosterman sums, 
and that they sometimes do not exist,
one is generally satisfied with  divisibility results and characterisation of Kloosterman sums modulo some integer (see \cite{MMeven, LisonekSeta, Div3, CHZ, Lisonek2009}). 

It is easy to see that binary Kloosterman sums are divisible by $4 = 2^2$, i.e.,
for all $a \in \f{2^n}$,
\begin{equation}\label{folklore}
\mathcal{K}_{2^n}(a) \equiv 0 \pmod{4}.
\end{equation}

They also satisfy (see \cite{LachWolf})
\[
-2^{n/2+1} \le \mathcal{K}_{2^n}(a) \le 2^{n/2+1},
\]
and take every value which is congruent to $0$ modulo $4$ in that range.

Helleseth and Zinoviev proved the following result which improved \eqref{folklore} 
one level higher, i.e., modulo $2^3$,  in the sense
of describing the $a$ for which $\mathcal{K}_{2^n}(a)$ is 0 or 4 modulo 8.

\bigskip

\begin{thm}\cite{HZ}\label{mod8}
For $a \in \f{2^n}$,
\[
\mathcal{K}_{2^n}(a) \equiv \left\{ \begin{array}{ll} 
0 \ (\mathrm{mod} \ 8) & \textrm{if } \Tr(a) = 0,\\
4 \ (\mathrm{mod} \ 8) & \textrm{if } \Tr(a) = 1. 
\end{array}
\right. 
\]
\end{thm}

Similar to the binary case, it is easy to see that ternary Kloosterman sums are divisible by $3$, i.e., for all $a \in \f{3^n}$,
\begin{equation}\label{folklore2}
\mathcal{K}_{3^n}(a) \equiv 0 \pmod{3}.
\end{equation}

Ternary Kloosterman sums satisfy (see Katz and Livn\'e  \cite{Katz})
\[
-2\sqrt{3^n} < \mathcal{K}_{3^n}(a) < 2\sqrt{3^n}
\]
and take every value which is congruent to $0$ modulo $3$ in that range.

We will prove the following theorem, a simple characterisation of ternary Kloosterman sums modulo $3^2$ using the trace map (similar to Helleseth-Zinoviev result for binary case), by using Stickelberger's theorem. 

\bigskip

\begin{thm}
For $a \in \f{3^n}$,
\[
\mathcal{K}_{3^n}(a) \equiv \left\{ \begin{array}{ll} 
0 \ (\mathrm{mod} \ 9) & \textrm{if } \Tr(a) = 0,\\
3 \ (\mathrm{mod} \ 9) & \textrm{if } \Tr(a) = 1,\\ 
6 \ (\mathrm{mod} \ 9) & \textrm{if } \Tr(a) = 2. 
\end{array}
\right. 
\]
\end{thm}

This result is implied by a result of van der Geer and van der Vlugt \cite{Geer}.

We will also give a characterisation modulo $3^3$ of Kloosterman sums, using the Gross-Koblitz formula. The characterisation will depend on a generalisation of the trace function. Note that the trace of an element $a \in \f{q}$ can be written as
\[
\Tr(a) := \sum_{i \in W_1} a^i,
\]

where $W_1 := \{ p^i \ | \ i \in \{0,\ldots,n-1\} \}$. We will use a generalised trace $\tau_S : \f{p^n} \to \f{p}$,
\[
\tau_S(a) := \sum_{i \in S} a^i,
\]
where $S$ can be assigned to any subset of $\{0,\ldots,p^n-2\}$ satisfying
\[
S^p := \{s^p \pmod{p^n-1} \ | \ s \in S\} = S,
\]
particularly quadratic and cubic powers of $p$, in contrast to the set of linear powers $W_1$. 

We will define the sets
\begin{align*}
X&:=\{r \in \{0,\dots,q-2\}|r = 3^i+3^j\},\ (i, j \text{ not necessarily distinct})\\
Y&:=\{r \in \{0,\dots,q-2\}|r = 3^i+3^j+3^k, i,j,k \text{ distinct}\},\\
Z&:=\{r \in \{0,\dots,q-2\}|r = 2\cdot3^i+3^j, i\ne j\}.
\end{align*}

Our main result is
\begin{thm}
Let $n\ge 3$, and let $q=3^n$. Then 
\begin{equation*}
\K_q(a)\equiv \left\{
	\begin{array}{rcccccl}
		0\pmod{27}\text{ if } &\Tr(a) = &0 &\text{ and } &\tau_{Y}(a) &+2\tau_{X}(a)&=0\\
		3\pmod{27}\text{ if } &\Tr(a) = &1 &\text{ and } &\tau_{Y}(a) & &= 2\\
		6\pmod{27}\text{ if } &\Tr(a) = &2 &\text{ and } &\tau_{Y}(a)&+\tau_{X}(a) &= 2\\
		9\pmod{27}\text{ if } &\Tr(a) = &0 &\text{ and } & \tau_{Y}(a)&+2\tau_{X}(a) &= 1\\
		12\pmod{27}\text{ if } &\Tr(a) = &1 &\text{ and } & \tau_{Y}(a)& &= 0\\
		15\pmod{27}\text{ if } &\Tr(a) = &2 &\text{ and } &\tau_{Y}(a)&+\tau_{X}(a) &= 0\\
		18\pmod{27}\text{ if } &\Tr(a) = &0 &\text{ and } &\tau_{Y}(a)&+2\tau_{X}(a) &=2\\
		21\pmod{27}\text{ if } &\Tr(a) = &1 &\text{ and } &\tau_{Y}(a)& &= 1\\
		24\pmod{27}\text{ if } &\Tr(a) = &2 &\text{ and } &\tau_{Y}(a)&+\tau_{X}(a) &= 1.
		\end{array} \right.
\end{equation*}

\end{thm}

Recently, we have proved a similar result for the binary case, using $\tau_{Q}$, where $Q := \{2^i + 2^j \ | \ i,j \in \{0,\ldots,n-1\}, i \ne j \}$.  
\bigskip

\begin{thm}\cite{arxivbinary}
For $a \in \f{2^n}$, 
\begin{equation*}
\mathcal{K}_{2^n}(a)\equiv \left\{
	\begin{array}{rccl}
			0 \pmod{16}& \text{ if } &\Tr(a) = 0 \text{ and } &\tau_{Q}(a) = 0,\\
			4\pmod{16}& \text{ if } &\Tr(a) = 1 \text{ and } &\tau_{Q}(a) = 1,\\
			8\pmod{16}& \text{ if } &\Tr(a) = 0 \text{ and } &\tau_{Q}(a) = 1,\\
			12\pmod{16}& \text{ if } &\Tr(a) = 1 \text{ and } &\tau_{Q}(a) = 0.
		\end{array} \right.
		\end{equation*}
\end{thm}

For the ternary case we mention a recent result due to Lisonek \cite{LisonekSeta} that gives a description of the elements $a \in \f{3^n}$ for which $\mathcal{K}(a) \equiv 0 \pmod{9}$, which is also implied by the van der Geer-van der Vlugt result.
\bigskip

\begin{thm}\cite{LisonekSeta}
Let $n \ge 2$. For any $a \in \f{3^n}, \mathcal{K}_{3^n}(a)$ is divisible by $9$ if and only if $\Tr(a) = 0$.  
\end{thm}

In Sections \ref{sect::Stickelberger} and \ref{sect::Fourier}, we will introduce the techniques we use. In Section \ref{results} we will give the results modulo $9$. In Section \ref{GKmod27} we will give the modulo $27$ result.

\section{Stickelberger's theorem}\label{sect::Stickelberger}

Let $p$ be a prime (in Section \ref{results} we set $p=3$).
Consider multiplicative characters taking their values in an algebraic extension of $\mathbb
{Q}_p$. Let $\xi$ be a primitive $(q-1)^{\text{th}}$ root of unity in a 
fixed algebraic closure of  $\mathbb{Q}_p$.  
The group of multiplicative characters of $\mathbb{F}_q$ (denoted $
\widehat{\mathbb{F}_q^{\times}}$) is cyclic of order $q-1$. 
The group $\widehat{\mathbb{F}_q^
{\times}}$ is generated by the Teichm\"uller character $\omega : \mathbb{F}_q \to 
\mathbb{Q}_p(\xi)$, which, for a fixed generator $t$ of $\mathbb{F}_q^{\times}$, is defined 
by $\omega (t^j) = \xi^j$. We set $\omega (0)$ to be $0$. An equivalent definition is that $\omega$ satisfies \[\omega(a) \equiv a \pmod{p}\] for all $a \in \f{q}$.

Let $\zeta$ be a fixed primitive $p$-th root of unity in the 
fixed algebraic closure of  $\mathbb{Q}_p$.  
Let $\mu$ be the canonical additive character of $\mathbb{F}_q$, \[\mu(x) = \zeta^{\Tr (x)}
\] 
where $\Tr$ denotes the absolute trace map from  $\mathbb{F}_q$ to  $\mathbb{F}_p$.

The Gauss sum (see \cite{LN,Wash}) of a character $\chi \in \widehat{\mathbb{F}_q^
{\times}}$  is defined as \[\tau(\chi) = -\sum_{x \in \mathbb{F}_q}\chi(x)\mu(x)\,.\]
We define \[g(j):=\tau(\omega^{-j})\,.\]

For any positive integer $j$, let $\wt_p(j)$ denote the $p$-weight of $j$, i.e.,
\[ \wt_p(j)=\sum_i j_i
\] 
where $\sum_i j_ip^i$ is the $p$-ary
expansion of $j$.

Let $\pi$ be the unique ($p-1$)th root of $-p$ in $\mathbb{Q}_p(\xi,\zeta)$ satisfying \[\pi \equiv \zeta - 1 \pmod{\pi^{2}}\,.\]
Wan \cite{Wan}
noted that the following improved version of Stickelberger's theorem is a direct consequence of the Gross-Koblitz formula (see Section \ref{GKmod27}).
\bigskip

\begin{thm}\cite{Wan}\label{Stick+}
Let $1 \le j < q-1$ and let $j =j_0+j_1p+\dots + j_{n-1}p^{n-1}$. Then 
\[g(j) \equiv \frac{\pi^{\wt_p(j)}}{j_{0}!\cdots j_{n-1}!} \pmod{\pi^{\wt_p(j)+p-1}}\,.\]
\end{thm}

Stickelberger's theorem, as usually stated,  is the same congruence modulo $\pi^{\wt_p(j)+1}$.

We have (see \cite{GK}) that
$(\pi)$ is the unique prime ideal of $\mathbb{Q}_{p}(\zeta, \xi)$ lying above $p$. 
Since $\mathbb{Q}_{p}(\zeta, \xi)$ is an unramified extension of $\mathbb{Q}_{p}(\zeta)$,
a totally ramified (degree $p-1$) extension of $\mathbb{Q}_{p}$, it follows that $(\pi)^{p-1} = (p)$
and $\nu_p(\pi) = \frac{1}{p-1}$.
Here $\nu_p$ denotes the $p$-adic valuation.

Therefore Theorem \ref{Stick+} implies that
$\nu_{\pi}(g(j)) = \wt_p(j)$, and because 
$\nu_{p}(g(j)) =\nu_{\pi}(g(j)) \cdot \nu_p (\pi)$ we get 
\begin{equation}\label{val}
\nu_p(g(j)) = \frac{\wt_p(j)}{p-1}\,.\end{equation}

In this paper we have $p=3$.  In that case, $\pi=-2\zeta-1$ and $\pi^2 = -3$. Hence \eqref{val} becomes
 \begin{equation}\label{val3}
\nu_3(g(j)) = \frac{\wt_3(j)}{2}\,.\end{equation}

\section{Fourier coefficients}\label{sect::Fourier}

The Fourier transform of a function $f:\mathbb{F}_q \to \mathbb{C}$ at $a \in \mathbb{F}_q$ 
is defined to be
\[\tfr{f}(a) = \sum_{x \in \mathbb{F}_q} f(x)\mu(ax)\,.\]
The complex number $\tfr{f}(a) $ is called the Fourier coefficient of $f$ at $a$.

Consider monomial functions defined by $f(x) = \mu(x^d)$. 
When $d=-1$ we have $\tfr{f}(a) = \mathcal{K}_{p^n}(a)$.
By a similar 
Fourier analysis argument to that in Katz \cite{KatzGaussSums} or
Langevin-Leander \cite{LL}, for any $d$ we have
\[
\tfr{f}(a) =\frac{q}{q-1}+ \frac{1}{q-1} \sum_{j=1}^{q-2}\tau(\bar{\omega}^j)\ \tau
(\omega^{jd})\ \bar{\omega}^{jd}(a) 
\]
and hence
\[
\tfr{f}(a) \equiv - \sum_{j=1}^{q-2}\tau(\bar{\omega}^j)\ \tau
(\omega^{jd})\ \bar{\omega}^{jd}(a) \pmod{q}\,.
\]
We will use this to obtain congruence information about  Kloosterman sums.
Putting $d=-1=p^n-2$, the previous congruence becomes

\begin{equation}\label{firstKL}
\mathcal{K}(a) \equiv - \sum_{j=1}^{q-2}(g(j))^2\ \omega^{j}(a) \pmod{q}.
\end{equation}
In this paper, $p=3$.
Equation (\ref{val3}) gives the $3$-adic valuation of the Gauss sums $g(j)$,
and the $3$-adic valuation of each term in  equation (\ref{firstKL}) follows.
Our proofs will consider (\ref{firstKL}) at various levels, i.e., modulo $3^2$ and $3^3$.

\section{Ternary Kloosterman sums modulo 9}\label{results}
 
In this section we will prove our result using Stickelberger's theorem. First we need a lemma which helps us in our proof.

\begin{lem}\label{SumInv}
Let $p$ be a prime, $q = p^n$ and $r \in \mathbb{F}_p^{\times}$. If $T_r$ denotes the set $\{a \in \mathbb{F}_q \mid \Tr (a) = r\}$, then \[\sum_{t \in T_r} t^{-1} = r^{-1}\,.\]
\end{lem}

\begin{proof}

Consider the polynomials \[g(x)= \prod_{t \in T_r}(x-t) \,,\]\[h(x)= \prod_{t \in 
T_r}(x-t^{-1}) \,.\]

Note that $g(x)$ vanishes on the $p^{n-1}$ elements of $T_r$. Thus 
\[
g(x) = x^{p^{n-1}}+x^{p^{n-2}}+\dots+x-r.
\]
In particular, $$\prod_{t \in T_r}(-t) = -r,$$ so $$\prod_{t \in T_r}(-t^{-1}) = -r^{-1}.$$

The reciprocal polynomial of $g$ is $g^{*}(x) = x^{p^{n-1}}g(1/x)$.

We therefore get
\begin{align*}
h(x) & = -r^{-1}g^{*}(x) \\
		 & = -r^{-1}x^{p^{n-1}}g(1/x)\\
     & = x^{p^{n-1}}-r^{-1}x^{p^{n-1}-1}-\dots-r^{-1}x^{p^{n-1}-p^{n-2}}-r^{-1}\,.
\end{align*} 
Thus \[\sum_{t \in T_r}(-t^{-1})=-r^{-1}\,.\]\qedhere

\end{proof}

From now on, we set $p = 3$, so that $\mathcal{K}_q (a)$ is an integer for $a \in \f{q}$. Since there will not be any confusion with binary Kloosterman sums we will write $\mathcal{K}(a)$ for $\mathcal{K}_q(a)$. We consider the function $f(x) = \mu(x^{-1}) = \mu(x^{q-2})$. Then $\tfr{f}(a)$ is the Kloosterman sum $\mathcal{K}(a)$. The following lemma will be needed.

\begin{lem}\label{cong3}
Let $q = 3^n$, and $T_1$ be as defined above. Then
\[\sum_{z\in T_1}\bar{\omega}(z) \equiv 1 \pmod{3}.\]
\end{lem}
\begin{proof}
Follows directly from Lemma \ref{SumInv} and the definition of the Teichm\"uller character.\qedhere
\end{proof}

We can now state our main result of this section.
\bigskip

\begin{thm}\label{mainthmfull}
Let $q=3^n$ for some integer $n>1$.
For $a \in \mathbb{F}_q$,
\begin{equation*}
\mathcal{K}_q(a)\equiv \left\{
	\begin{array}{ll}
			0 \pmod{9} & \text{if } \Tr (a) = 0,\\
			3 \pmod{9} & \text{if } \Tr (a) = 1,\\
			6 \pmod{9} & \text{if } \Tr (a) = 2.
		\end{array} \right. 
\end{equation*}
\end{thm}

\begin{proof}
 By \eqref{firstKL}
\begin{equation}
\mathcal{K}(a) \equiv - \sum_{j=1}^{q-2}g(j)^2\ \omega^{j}(a)\pmod{q}\,.\label{Kmodq}
\end{equation}

Let, for any $0 < t < q-1$, the $3$-adic expansion of $t$ be $t = t_{0}+3t_{1}+\dots+3^{n-1} t_{n-1}$ and let $\mathcal{P}$ be the prime of $\mathbb{Q}_3(\xi,\zeta)$ lying above $3$. 
As we mentioned in Section 2, Stickelberger's theorem implies that 
\begin{eqnarray}
\nu_{\mathcal{P}}(g(t)) & = & \wt_3(t)=t_{0}+t_{1}+\dots+t_{n-1} \nonumber\\
\nu_3(g(t)) & = & \frac{\wt_3(t)}{2}, \nonumber\\
\textrm{and so}\ \nu_3((g(t))^2) & = & \wt_3(t).\label{wt}
\end{eqnarray}
Now (\ref{wt}) implies that any term in the sum in \eqref{Kmodq} with
$\wt_3(j)>1$ will be 0 modulo 9, so \eqref{Kmodq} modulo 9 becomes
a sum over terms of weight 1 only:
\begin{equation*}
\mathcal{K}(a) \equiv - \sum_{0\leq i<n}g(3^i)^2\ \omega^{3^i}(a) \pmod{9}\,.
\end{equation*}
By Lemma 6.5 of \cite{Wash}, $g(3^i) = g(1)$, so we obtain
\begin{equation}\label{ka1}
\mathcal{K}(a) \equiv - g(1)^2\ \sum_{0\leq i<n}\omega^{3^i}(a)\pmod{9}\,.
\end{equation}
By  definition of $\omega$, we have 
\begin{equation}\label{ka2}
\sum_{0\leq i<n}\omega^{3^i}(a) \equiv \Tr (a) \pmod{3}\,.
\end{equation}
Since $\nu_3(g(1)^2) =\wt_3(1)=1$, 
the proof of the theorem reduces to determining $g(1)^2\  \bmod{9}$. 
We calculate, using the notation of Lemma \ref{SumInv},
\begin{align*}
g(1)& = -\sum_{x\in \mathbb{F}_q^{\times}}\bar{\omega}(x)\zeta^{\Tr (x)}\\
&=-\sum_{x\in T_0}\bar{\omega}(x) -\sum_{x\in T_1}\bar{\omega}(x)\zeta -\sum_{x\in 
T_1}\bar{\omega}(-x)\zeta^2\\
&=(\zeta^2-\zeta)\sum_{x\in T_1}\bar{\omega}(x)
\end{align*}
because $\bar{\omega}(-x)=-\bar{\omega}(x)$, $T_2=-T_1$,  and
the sum over $T_0$ is 0.
This implies
\[
 g(1)^2 = (\zeta^2-\zeta)^2\left(\sum_{x\in T_1}\bar{\omega}(x)\right)^2.
\]
But we have $(\zeta^2-\zeta)^2 = -3$. This, together with Lemma \ref{cong3}, implies 
\begin{equation}\label{6mod9}
g(1)^2 \equiv 6 \pmod{9}.
\end{equation}
Combining this with (\ref{ka2}), the congruence (\ref{ka1}) becomes
\begin{equation*}
\mathcal{K}(a) \equiv 3\ \Tr (a)\pmod{9}
\end{equation*}
as required.
  \qedhere
\end{proof}

Garaschuk and Lisonek proves the following theorem which characterises ternary Kloosterman sums modulo $2$.

\begin{thm}\label{mod2} \cite{Lisonek2009}
Let $\sqrt{a}$ denote any $b \in \f{3^n}$ such that $b^2=a$.  
\[
\mathcal{K}_{3^n}(a) \equiv \left\{ \begin{array}{ll} 
0 \ (\mathrm{mod} \ 2) & \textrm{if } a = 0 \textrm{ or } a \textrm{ is a square and } \Tr(\sqrt{a}) \ne 0,\\
1 \ (\mathrm{mod} \ 2) & \textrm{otherwise.}
\end{array}
\right. 
\]
\end{thm}

Theorem \ref{mainthmfull} and Theorem \ref{mod2} together give a full characterisation of ternary Kloosterman sums modulo $18$, which we summarise in the following corollary.

\begin{cor}\label{mod18}
Let $q=3^n$. For $a \in \mathbb{F}_q^\times$,
\begin{equation*}
\mathcal{K}_q(a)\equiv \left\{
	\begin{array}{rcccccl}
			0  \pmod{18} &\textit{ if } &\Tr(a) = 0 & \textit{and }  a & \textsf{square}      & \textit{with} &\Tr(\sqrt{a}) \ne 0,\\
			3  \pmod{18} &\textit{ if } &\Tr(a) = 1 & \textit{and }  a & \textsf{non-square}  & \textit{or}   &\Tr(\sqrt{a}) = 0, \\
			6  \pmod{18} &\textit{ if } &\Tr(a) = 2 & \textit{and }  a & \textsf{square}      & \textit{with} &\Tr(\sqrt{a}) \ne 0,\\
			9  \pmod{18} &\textit{ if } &\Tr(a) = 0 & \textit{and }  a & \textsf{non-square}  & \textit{or}   &\Tr(\sqrt{a}) = 0,  \\
			12 \pmod{18} &\textit{ if } &\Tr(a) = 1 & \textit{and }  a & \textsf{square}      & \textit{with} &\Tr(\sqrt{a}) \ne 0,\\
			15 \pmod{18} &\textit{ if } &\Tr(a) = 2 & \textit{and }  a & \textsf{non-square}  & \textit{or}   &\Tr(\sqrt{a}) = 0.
		\end{array} \right. 
\end{equation*}
\end{cor}

\section{Ternary Kloosterman sums modulo 27}\label{GKmod27}

To be able to give higher level congruences we will need a result stronger than Stickelberger's 
theorem.   Recall that Gauss sums lie in $\mathbb{Z}_{p}[\zeta, \xi]$,
and that $(\pi)$ is the unique prime ideal of $\mathbb{Z}_{p}[\zeta, \xi]$ lying above $p$. 
All congruences involving Gauss sums take place in this ring, so when we write
$g(j)^2 \equiv 6 \pmod{27}$ we mean that $g(j)^2-6$ is in the ideal $(27)$.
The {\bf Gross-Koblitz formula} \cite{GK,Rob} states that 
\begin{equation}\label{GKformula}
g(j) = \pi^{\wt_p(j)}\prod_{i=0}^{n-1}\Gamma_{p}\left(\left\langle \frac{p^{i}j}{q-1} \right\rangle\right)
\end{equation}
where $\langle x \rangle$ is the fractional part of a rational number $x$, and $\Gamma_p$ is the $p$-adic Gamma function $\Gamma_{p}:\mathbb{N} \to \mathbb{N}$ defined by (cf. \cite{Morita})

\[
\Gamma_{p}(k) = (-1)^k\prod_{\substack{t<k\\(t,p) = 1}}t \,.
\]

The following result helps one computing the $p$-adic Gamma function modulo $p^k$.

\begin{thm}[Generalised Wilson's theorem]\cite{Morita}

Suppose $x \equiv y \pmod{p^k}$. If $p^k \ne 4$, then
\[
\Gamma_{p}(x) \equiv \Gamma_p(y) \pmod{p^k}.
\]
\end{thm}

This theorem is actually a consequence of Gauss' generalisation of Wilson's theorem. Now let us prove a lemma on evaluations of the $p$-adic Gamma function. This lemma will allow us to 
evaluate Gauss sums for higher moduli and find Kloosterman congruences modulo $27$.
\bigskip

\begin{lem}\label{gammamod}Let $q=3^n$ and let $i$ be an integer in the range $[0,n-1]$. Then
\[
\Gamma_3\left(\left\langle \frac{3^{i}}{q-1} \right\rangle\right) \equiv 
\left\{ \begin{array}{ll}
13 \pmod{27} & \textrm{if } i = 1, \\
1 \pmod{27} & \textrm{if } i > 1.
\end{array} \right.
\]
\end{lem}

\begin{proof}
For any $j$, we have $3^j \le q$, and
\[
\left\langle\frac{3^i}{q-1}\right\rangle = \frac{3^i}{q-1} \equiv 3^i(3^j-1) \pmod{3^j},\] so 
\[
\Gamma_3\left(\left\langle \frac{3^{i}}{q-1} \right\rangle\right) \equiv \Gamma_3(26\cdot3^i) \pmod{27}.\]
If $i \ge 3$, then $26\cdot 3^i \equiv 0 \pmod{27}$, and
\[ 
\Gamma_3\left(\left\langle \frac{3^{i}}{q-1} \right\rangle\right) \equiv 1 \pmod{27}\,,
\]

Now $\Gamma_3(26\cdot3) \equiv \Gamma_3(24) \pmod{27}$ using Generalised Wilson's theorem. And $\Gamma_3(24) \equiv 13 \pmod{9}$. Similarly:
\begin{eqnarray*}
\Gamma_3(26\cdot9) & \equiv & 1 \pmod{27}.
\end{eqnarray*}
\end{proof}

Lemma \ref{gammamod} allows us to compute Gauss sums modulo $27$:
\bigskip

\begin{lem}\label{wt1lem}
Let $q=3^n$. Then
\[
g(j)^2 \equiv \left\{ \begin{array}{ll}
6 \pmod{27} & \textrm{if } \wt_p(j) =  1, \\
9 \pmod{27} & \textrm{if } \wt_p(j) =  2, \\
0 \pmod{27} & \textrm{if } \wt_p(j) \ge 3.
\end{array} \right.
\]
\end{lem}

\begin{proof}
Suppose $\wt_p(j) = 1$. By the Gross-Koblitz formula and Lemma \ref{gammamod}, 
\[
g(j) \equiv 13 \pi \pmod{27}.
\]
Let 
\[
g(j) = 27A + 13\pi
\]
for some $A\in \mathbb{Z}_{p}[\zeta, \xi]$. Then
\begin{align*}
g(j)^2 & = 27^2 A^2 + 2\cdot 27\cdot 13 A  + 169\pi^{2} \\
			 & \equiv 169\pi^{2} \pmod{27} \\
			 & \equiv 6 \pmod{27}
\end{align*}
since $\pi^2=-3$. Now suppose $\wt_p(j) = 2$. By the Gross-Koblitz formula, 
\[
g(j) \equiv -3 \pmod{9}. 
\] Thus $g(j) = 9X-3$ for some $X\in \mathbb{Z}_{p}[\zeta, \xi]$,
so \[g(j)^2 = 81X^2-54X+9 \equiv 9 \pmod{27}.\]

It is clear from the Gross-Koblitz formula that if $\wt_p(j)>2$, then
\[
27 | \pi^{2\wt_p(j)} | g(j)^2.
\]
\end{proof}

Consider again the trace function $\Tr : \f{p^n} \to \f{p}$,
\[
\Tr (c) = c + c^p + c^{p^2} + \dots + c^{p^{n-1}}.
\]

We wish to generalise this definition to a larger class of finite field sums, which includes the usual trace function as a special case.
\bigskip

\begin{defn}\label{trace}
Let $p$ be a prime, let $n\ge 1$ be an integer and let $q=p^n$.
For any $S \subseteq \mathbb{Z}/(q-1)\mathbb{Z}$ satisfying $S^p=S$
where
$S^p := \{s^p \ | \ s \in S\}$,
define the $S$-\emph{trace} to be the function $\tau_S:\f{q} \to \f{p}$,
\[\tau_S(c) := \sum_{s\in S} c^s\,.\]
\end{defn}

Let 
\begin{align*}
X&:=\{r \in \{0,\dots,q-2\}|r = 3^i+3^j\},\ (i, j \text{ not necessarily distinct})\\
Y&:=\{r \in \{0,\dots,q-2\}|r = 3^i+3^j+3^k, i,j,k \text{ distinct}\},\\
Z&:=\{r \in \{0,\dots,q-2\}|r = 2\cdot3^i+3^j, i\ne j\}.
\end{align*}

Now we are ready to prove our result on Kloosterman sums modulo $27$.
\bigskip

\begin{thm}\label{mod27} Let $\K_q$ be the usual $q$-ary Kloosterman sum, 
let
\[\tfr{\Tr}(a)=\sum_{\wt_{3}(i) = 1}\omega^{i}(a), \ \textrm{and let}\ \
\tfr{\tau_{X}}(a) = \sum_{\wt_{3}(j) = 2}\omega^j(a).\] Then

\begin{equation}\label{hateq}
\K_{3^n}(a) \equiv 21\tfr{\Tr}(a)+18\tfr{\tau_{X}}(a) \pmod{27}.
\end{equation}
\end{thm}
\begin{proof}
Using   \eqref{firstKL} and Lemma \ref{wt1lem}, we get
\begin{align*}
\mathcal{K}(a) & \equiv - \sum_{j=1}^{q-2}g(j)^2\ \omega^{j}(a)\pmod{q} \\
		& \equiv - \sum_{\wt_{3}(j)=1}g(j)^2 {\omega}^{j}(a) - \sum_{\wt_{3}(j)=2} g(j)^2 {\omega}^{j}(a) \pmod{27} \\
		& \equiv -6\sum_{\wt_{3}(j)=1}{\omega}^{j}(a) - 9\sum_{\wt_{3}(j)=2}{\omega}^{j}(a) \pmod{27} \\
&\equiv 21\tfr{\Tr}(a)+18\tfr{\tau_{X}}(a) \pmod{27}.\qedhere
 \end{align*}
\end{proof}

It would be preferable to express the above result in terms of operations within $\f{q}$ itself. Note that in (\ref{hateq}) we only need $\tfr{\Tr}(a)$ modulo $9$ and $\tfr{\tau_{X}}(a)$ modulo $3$. We have
\[
\tau_{X} (a) \equiv \tfr{\tau_{X}}(a) \pmod{3}.
\]
We need to find some condition for $\tfr{\Tr}(a)$ modulo $9$ using functions from $\f{q}$ to $\f{p}$. We will do that in the proof of the following corollary.
\bigskip

\begin{cor}\label{fieldsums}
Let $n\ge 3$, and let $q=3^n$. Then 
\[\K_q (a) \equiv 21 \Tr(a)^3+18\tau_{Z}(a)+9\tau_{Y}(a)+18\tau_{X}(a) \pmod{27}.\]
\end{cor}
\begin{proof}

First note that $\tfr{Q}(a) \equiv \tau_{X}(a) \pmod{3}$, by the basic property of the Teichm\"uller character.

To determine $\tfr{\Tr}(a) \bmod{9}$, we compute 
\begin{align*}
\tfr{\Tr}(a)^3 &= \sum_{i,j,k \in \{0,\dots,n-1\}}\omega(a^{3^i+3^j+3^k})\\
&=\tfr{\Tr}(a)+3\sum_{r\in Z}\omega(a^{r})+6\sum_{r \in Y}\omega(a^{r})\,,
\end{align*}
and note the elementary fact that if $x \equiv y \pmod{m}$, then $x^m \equiv y^m \pmod{m^2}$. This means that $\tfr{\Tr}(a)^3 \bmod{9}$ is given by $\tfr{\Tr}(a) \bmod{3} = \Tr (a)$, i.e. $\tfr{\Tr}(a)^3 \bmod{9} = \Tr(a)^3$.

Since \[\sum_{r\in Z}\omega(a^{r}) \equiv \tau_{Z}(a)\pmod{3}\]and \[\sum_{r \in Y}\omega(a^{r}) \equiv \tau_{Y}(a)\pmod{3}\,,\]we have that

\[\tfr{\Tr}(a) \equiv \Tr(a)^3 - 3\tau_{Z}(a) - 6\tau_{Y}(a)\pmod{9},\]proving the result.
\end{proof}

Note that
\[\Tr(a) \tau_{X}(a) = \Tr(a)+2\tau_{Z}(a)\,.\]
Thus Corollary \ref{fieldsums} can be rewritten as
\begin{equation}\label{final27}
\K_q (a) \equiv 21 \Tr(a)^3+18\Tr(a)+18\tau_{X}(a)+9\Tr(a)\tau_{X}(a)+9\tau_{Y}(a) \pmod{27}.
\end{equation}
The smallest field for which each of the 27 possible values of $(\Tr(a), \tau_{X}(a), \tau_{Y}(a))$ occurs is $\f{3^6}$.
\bigskip

\begin{cor}
Let $n\ge 3$, and let $q=3^n$. Then 
\begin{equation*}
\K_q(a)\equiv \left\{
	\begin{array}{rcccccl}
		0\pmod{27}\text{ if } &\Tr(a) = &0 &\text{ and } &\tau_{Y}(a) &+2\tau_{X}(a)&=0\\
		3\pmod{27}\text{ if } &\Tr(a) = &1 &\text{ and } &\tau_{Y}(a) & &= 2\\
		6\pmod{27}\text{ if } &\Tr(a) = &2 &\text{ and } &\tau_{Y}(a)&+\tau_{X}(a) &= 2\\
		9\pmod{27}\text{ if } &\Tr(a) = &0 &\text{ and } & \tau_{Y}(a)&+2\tau_{X}(a) &= 1\\
		12\pmod{27}\text{ if } &\Tr(a) = &1 &\text{ and } & \tau_{Y}(a)& &= 0\\
		15\pmod{27}\text{ if } &\Tr(a) = &2 &\text{ and } &\tau_{Y}(a)&+\tau_{X}(a) &= 0\\
		18\pmod{27}\text{ if } &\Tr(a) = &0 &\text{ and } &\tau_{Y}(a)&+2\tau_{X}(a) &=2\\
		21\pmod{27}\text{ if } &\Tr(a) = &1 &\text{ and } &\tau_{Y}(a)& &= 1\\
		24\pmod{27}\text{ if } &\Tr(a) = &2 &\text{ and } &\tau_{Y}(a)&+\tau_{X}(a) &= 1.
		\end{array} \right.
\end{equation*}

\end{cor}

\begin{proof}
Restatement of equation \ref{final27}.
\end{proof}

The Kloosterman sums modulo $54$ can be given by combining (\ref{fieldsums}) and Theorem \ref{mod2}.

\bibliographystyle{plain}
\bibliography{main}
\end{document}